\documentclass[11pt]{amsart}
\usepackage{mathtools}
\usepackage{amssymb}
\usepackage{amsmath}
\usepackage{mathrsfs}
\usepackage[pagebackref,colorlinks=true,linkcolor=blue,citecolor=blue]{hyperref}
\usepackage{etoolbox}
\usepackage{enumitem}
\usepackage{tikz-cd}
\usepackage[normalem]{ulem}
\usepackage{color}

\UseRawInputEncoding

\newtheorem{theorem}{Theorem}[section]
\newtheorem{lemma}[theorem]{Lemma}
\newtheorem{proposition}[theorem]{Proposition}

\newtheorem{conjecture}[theorem]{Conjecture}
\newtheorem{definition}[theorem]{Definition}

\newtheorem{desideratum}[theorem]{Desideratum}

\newtheorem{remark}[theorem]{Remark}

\newcommand{\BC}{\mathbb{C}}

\newcommand{\ol}{\overline}

\newcommand{\inv}{^{-1}}

\newcommand{\CO}{\mathcal{O}}

\newcommand{\mfr}[1]{\mathfrak{#1}}

\newcommand{\RG}{\mathrm{G}}
\newcommand{\SL}{\mathrm{SL}}

\newcommand{\SO}{\mathrm{SO}}
\newcommand{\Sp}{\mathrm{Sp}}
\newcommand{\RU}{\mathrm{U}}

\newcommand{\Ind}{\mathrm{Ind}}
\newcommand{\Jac}{\mathrm{Jac}}
\newcommand{\LLC}{\mathrm{LLC}}
\newcommand{\Ka}{\operatorname{Ka}}
\newcommand{\AZ}{\operatorname{AZ}}

\DeclareMathOperator{\Ad}{Ad}
\DeclareMathOperator{\Hom}{Hom}

\DeclareMathOperator{\im}{im}

\numberwithin{equation}{section}

\let\oldbullet\bullet
\renewcommand{\bullet}{{\vcenter{\hbox{\tiny$\oldbullet$}}}}

\title[Generic case of Jiang's conjecture]{On the generic case of Jiang's conjecture}

\author[C.-H. Lo]{Chi-Heng Lo}
\address{Department of Mathematics\\
National University of Singapore\\
119076, Singapore}
\email{{ch\_lo@nus.edu.sg}}

\subjclass[2020]{Primary 11F70, 22E50; Secondary 11F85, 22E55}

\date{\today}

\keywords{Admissible Representations, Aubert-Zelevinsky Involution, Generic Representations, Local Langlands Correspondence, Wavefront Sets}

\begin{document}

\begin{abstract}
    Jiang's conjecture predicts that if $\pi$ is a generic representation of a quasi-split $p$-adic group, then the $L$-parameter of the Aubert-Zelevinsky involution of $\pi$ must be trivial on $\SL_2$. In this paper, we prove this statement for a general quasi-split group under certain assumptions on the conjectural local Langlands correspondence.
\end{abstract}

\maketitle

%\tableofcontents

\section{Introduction}\label{sec intro}

Let $F$ be a non-Archimedean local field of characteristic zero, let $\RG$ be a connected reductive algebraic group defined over $F$, and put $G= \RG(F)$. Given an irreducible admissible representation $\pi$ of $G$, the Harish-Chandra-Howe local character expansion (\cite{HC78}) expresses the character $\Theta(\pi)$, on a sufficiently small neighborhood of the identity, as a linear combination
\[ \Theta(\pi)(\exp(X)) = \sum_{\CO} c_{\CO}(\pi)\, \widehat{\mu}_{\CO}(X)\]
of Fourier transforms of nilpotent orbital integrals, the sum running over the nilpotent orbits $\CO$ in $\mathfrak{g}(F)$. Let $\mathfrak{n}^{m}(\pi)$ be the set of maximal $\CO$ (with respect to the topological closure ordering) with $c_{\CO}(\pi)\neq 0$, and let $\textrm{WF}^{geo}(\pi)$ be the corresponding set of geometric nilpotent orbits over $\ol{F}$, taken maximal under the closure ordering. The set $\textrm{WF}^{geo}(\pi)$ is called the \emph{(geometric) wavefront set} of $\pi$.

By the work of M{\oe}glin-Waldspurger (\cite{MW87}, for residual characteristic $\neq 2$) and Varma (\cite{Var14}, for residual characteristic $2$), the maximal orbits $\CO$ for which the coefficient $c_{\CO}(\pi)\neq 0$ 
coincide with the maximal nilpotent orbits admitting a nonzero generalized Whittaker model of \(\pi\).
 As a special case, taking $\CO$ regular, $\pi$ is generic with respect to some Whittaker datum if and only if its geometric wavefront set is the regular orbit, that is,
\begin{align}\label{eq rodier intro}
    \pi \text{ is generic } \Longleftrightarrow \textrm{WF}^{geo}(\pi) = \{\CO_{reg}\}.
\end{align}

\subsection{Jiang's conjecture and the upper bound conjecture}\label{subsec intro jiang}

The Shahidi conjecture says that a tempered $L$-packet of a quasi-split group contains a generic member. Jiang (\cite[Conjecture 4.2]{Jia14}; see also \cite[Conjecture 1.7]{LS23}) proposed a generalization to arbitrary local Arthur packets, which predicts an upper bound for the wavefront sets of their members: for any local Arthur parameter $\psi$ of $G$, any $\pi \in \Pi_{\psi}$ and any $\CO \in \textrm{WF}^{geo}(\pi)$,
\[ \CO \leq d_{BV}(\CO_{\psi}),\]
where $\CO_{\psi}$ is the geometric nilpotent orbit attached to the restriction of $\psi$ to the Arthur $\SL_2$, and $d_{BV}$ is the Barbasch-Vogan duality map (\cite{BV85}).

In \cite{HLLS24}, and independently in \cite{CK25}, this prediction was restated in a form that does not refer to Arthur packets and that applies to \emph{every} irreducible admissible representation. To an $L$-parameter $\phi$ of $G$ one attaches a nilpotent orbit $\CO_{\phi}$ in $\widehat{\mathfrak{g}}(\BC)$, namely the orbit of the nilpotent element in an $\mathfrak{sl}_2$-triple associated to $\phi|_{\SL_2(\BC)}$.

\begin{conjecture}[Upper bound conjecture, {\cite[Conjecture 1.1]{HLLS24}}, \cite{CK25}]\label{conj upper bound}
    Assume the local Langlands correspondence for $G$. Let $\pi$ be an irreducible admissible representation of $G$. Then for any $\CO \in \textrm{WF}^{geo}(\pi)$,
    \[ \CO \leq d_{BV}\left(\CO_{\phi_{\widehat{\pi}}}\right),\]
    where $\widehat{\pi}$ is the Aubert-Zelevinsky involution of $\pi$ and $\phi_{\widehat{\pi}}$ its $L$-parameter.
\end{conjecture}

Under the assumptions on the theories of local Arthur packets and local ABV-packets stated in \cite{HLLS24}, Conjecture \ref{conj upper bound} is equivalent to Jiang's conjecture for representations of Arthur type, and to the analogous statement for local ABV-packets; see \cite[Theorems 1.4, 1.6]{HLLS24}. Ciubotaru and Kim verified it in a number of cases, including the depth-zero supercuspidal representations of classical groups and all irreducible representations of $\mathrm{G}_2$ (\cite{CK25}). More recently, Atobe and Ciubotaru (\cite{AC26}) established Conjecture \ref{conj upper bound} for all irreducible admissible representations of split $\Sp_{2n}$ and split $\SO_{2n+1}$ when $p \gg 0$, using Waldspurger's work on endoscopic transfer.

\subsection{The generic case}\label{subsec intro generic}

This paper is concerned with the case of Conjecture \ref{conj upper bound} in which $\pi$ is generic. We recall the discussion in \cite[Remark 1.5]{CK25}. Combining \eqref{eq rodier intro} with Conjecture \ref{conj upper bound} forces
\[ \CO_{reg}\leq d_{BV}\left(\CO_{\phi_{\widehat{\pi}}}\right),\]
and hence $d_{BV}(\CO_{\phi_{\widehat{\pi}}}) = \CO_{reg}$, since the regular orbit is maximal. This holds precisely when $\CO_{\phi_{\widehat{\pi}}}$ is the zero orbit, that is, when $\phi_{\widehat{\pi}}$ is trivial on $\SL_2(\BC)$. In other words, the generic case of Jiang's conjecture is equivalent to the following statement, which is the subject of this paper. Here $\Pi(G)$ denotes the set of equivalence classes of irreducible admissible representations of $G$, and $\Phi(G)$ the set of equivalence classes of $L$-parameters of $G$.

\begin{conjecture}\label{conj main}
    Let $G$ be quasi-split and equipped with a local Langlands correspondence
    \[ \LLC: \Pi(G) \longrightarrow \Phi(G).\]
    If $\pi \in \Pi(G)$ is generic, then $\LLC(\widehat{\pi})$ is trivial on $\SL_2(\BC)$.
\end{conjecture}

So Conjecture \ref{conj main} is a necessary condition for Conjecture \ref{conj upper bound}. For the quasi-split classical groups $\RU_n, \Sp_{2n}, \SO_{2n+1}, \SO_{2n}^{\alpha}$ and their pure inner forms, Conjecture \ref{conj main} was proved in \cite{LLS24}, using the explicit description of the local Langlands correspondence for these groups, especially the explicit parametrization of generic supercuspidal representations. In this paper we prove Conjecture \ref{conj main} for a general quasi-split group, by assuming properties of the conjectural local Langlands correspondence instead of restricting the group.

Our main result is the following. The desiderata \textbf{(LC)} and \textbf{(IT)} are stated in Desideratum \ref{des LC IT}. They say that $\LLC$ is compatible with the Langlands classification and with tempered parabolic induction, and both are expected to hold for any local Langlands correspondence.

\begin{theorem}\label{thm intro main}
    Let $G$ be quasi-split. Assume that there exists a local Langlands correspondence $\LLC_{\mathcal{L}(G)}$ for $G$ and its standard Levi subgroups  satisfying \emph{\textbf{(LC)}} and \emph{\textbf{(IT)}}. Assume in addition that one of the following holds.
    \begin{enumerate}
        \item [(a)] $p \gg 0$ with respect to $G$, and $\LLC_{M}(\pi)= \LLC^{\Ka}_{M}(\pi)$ for every $M \in \mathcal{L}(G)$ and every non-singular supercuspidal representation $\pi$ of $M$ in the sense of \cite{Kal21}, where $\LLC^{\Ka}$ is Kaletha's correspondence of \cite{Kal21}.
        \item [(b)] $\LLC_{\mathcal{L}(G)}$ satisfies the Gross-Prasad-Rallis conjecture (Conjecture \ref{conj GPR}), the uniqueness of the generic member in a tempered $L$-packet (Conjecture \ref{conj uniqueness}), and Haines' conjecture (Conjecture \ref{conj Haines}).
    \end{enumerate}
    Then Conjecture \ref{conj main} holds.
\end{theorem}

\subsection{Comments on the proof}\label{subsec intro proof}

The proof splits into a reduction step and a supercuspidal step, carried out in \S \ref{sec reduction} and \S \ref{sec supercuspidal} respectively.

The reduction is the following theorem, proved in \S \ref{sec reduction}.

\begin{theorem}\label{thm reduction}
    Assume that there exists a local Langlands correspondence $\LLC_{\mathcal{L}(G)}$ for $G$ and its standard Levi subgroups  satisfying \emph{\textbf{(LC)}} and \emph{\textbf{(IT)}}. If Conjecture \ref{conj main} holds for all generic supercuspidal representations of Levi subgroups of $G$, then it holds for all generic representations of $G$.
\end{theorem}

The main ingredients of the proof is 
the compatibility between the Jacquet module and the Aubert-Zelevinsky involution (Proposition \ref{prop AZ jacquet}) and  the generalized injectivity conjecture of Casselman-Shahidi (\cite[Conjecture 1.2]{CS98}). The latter is established through a sequence of works of Casselman-Shahidi, Hanzer, and Dijols, and finally proved in full generality by Solleveld; see Remark \ref{rmk injectivity} for details.

The supercuspidal case is where the extra hypotheses are needed, and we give two independent arguments, one for each of the conditions (a) and (b) of Theorem \ref{thm intro main}. The first argument (\S \ref{subsec kaletha}), which assumes condition (a), is short: by \cite{KLT}, for $p \gg 0$ every generic supercuspidal representation is non-singular in the sense of \cite{Kal21}, and Kaletha's construction attaches to such a representation a discrete parameter that is trivial on $\SL_2(\BC)$. The second argument (\S \ref{subsec GPR Haines}), which assumes condition (b), has no assumptions on $p$ and proceeds by contradiction. Let $\pi$ be a generic supercuspidal representation and write  $\phi = \LLC(\pi)$. If $\phi$ were not trivial on $\SL_2(\BC)$, the Gross-Prasad-Rallis and Haines conjectures would let us construct a non-supercuspidal generic member of the $L$-packet $\Pi_{\phi}$. Since $\pi$ is another generic member of $\Pi_{\phi}$, this contradicts the uniqueness of the generic member in a tempered $L$-packet. We conclude that $\phi$ must be trivial on $\SL_2(\BC)$.

\subsection{Organization}\label{subsec intro organization}

In \S \ref{sec prelim} we fix notation and recall what we need about Whittaker data and generic representations, the Langlands classification and the standard module conjecture, the Aubert-Zelevinsky involution, and the local Langlands correspondence. In \S \ref{sec reduction} we prove the reduction to the supercuspidal case. In \S \ref{sec supercuspidal} we treat the supercuspidal case under each of the two sets of hypotheses, which finishes the proof of Theorem \ref{thm intro main}.

\subsection*{Acknowledgements}
The author would like to thank Baiying Liu and Freydoon Shahidi for their constant support, and Cheng-Chiang Tsai for helpful discussions. 

\subsection*{Use of AI}
In preparing this paper, the author used Claude to transcribe the handwritten manuscript into \LaTeX, to draft the introduction (\S \ref{sec intro}) and preliminaries (\S \ref{sec prelim}), and for editing throughout.

\section{Preliminaries}\label{sec prelim}

Throughout this paper, $F$ denotes a non-Archimedean local field of characteristic zero, with Weil group $W_F$ and normalized absolute value $|\cdot|$. We let $\RG$ denote a connected reductive algebraic group defined and quasi-split over $F$, and put $G=\RG(F)$. We fix a Borel subgroup $B=TU$ of $\RG$ defined over $F$ and let $\mathcal{L}(G)$ denote the collection of standard Levi subgroups of $G$. We let $\Pi(G)$ denote the set of equivalence classes of irreducible smooth representations of $G$, and $\Phi(G)$ the set of equivalence classes of $L$-parameters of $G$, that is, of admissible homomorphisms
\[ \phi: W_F \times \SL_2(\BC) \longrightarrow {}^{L}G = \widehat{G}(\BC)\rtimes W_F\]
taken up to $\widehat{\RG}(\BC)$-conjugacy. For a Levi subgroup $M$ of $G$, the natural embedding ${}^{L}M \hookrightarrow {}^{L}G$ is well defined up to conjugacy. We write $\pi \leq \tau$ whenever $\pi$ is a subquotient of $\tau$.

\subsection{Whittaker data and generic representations}\label{subsec whittaker}

Let $\mathcal{U}=U(F)$ be the unipotent radical of $B(F)$. A character $\psi$ of $\mathcal{U}$ is called \emph{non-degenerate} (or generic) if its restriction to the root subgroup attached to each simple root is non-trivial. We fix a Whittaker datum $\mfr{w}$ of $G$ defined as follows throughout this paper.

\begin{definition}\label{def whittaker datum}
    A \emph{Whittaker datum} of $G$ is a $G$-conjugacy class $\mfr{w}= (B,\psi)$ of pairs consisting of a Borel subgroup $B$ of $\RG$ defined over $F$ and a non-degenerate character $\psi$ of $\mathcal{U}$. An irreducible admissible representation $\pi \in \Pi(G)$ is called \emph{$\mfr{w}$-generic} if
    \[ \Hom_{\mathcal{U}}(\pi, \psi)\neq 0,\]
    and \emph{generic} if it is $\mfr{w}$-generic for at least one Whittaker datum $\mfr{w}$.
\end{definition}

The fixed Whittaker datum $\mfr{w}$ of $G$ determines a Whittaker datum $\mfr{w}_M$ for each $M \in \mathcal{L}(G)$, characterized by the following property: for a parabolic subgroup $P$ with Levi subgroup $M$ and $\sigma \in \Pi(M)$, the induced representation $\Ind_{P}^{G}\sigma$ is $\mfr{w}$-generic if and only if $\sigma$ is $\mfr{w}_M$-generic. This is the heredity of Whittaker models, due to Rodier (\cite{Rod73}; see also \cite[Corollary 1.7]{CS80} and \cite[\S 3]{Sha90}). We shall use the following refinement.

\begin{proposition}[\cite{Rod73}]\label{prop rodier}
    Let $P=MN$ be a standard parabolic subgroup of $G$ and let $\sigma \in \Pi(M)$ be $\mfr{w}_M$-generic. Then $\Ind_{P}^{G}\sigma$ admits a unique $\mfr{w}$-generic irreducible subquotient, and it occurs with multiplicity one in $\Ind_{P}^{G}\sigma$.
\end{proposition}

\subsection{Langlands classification and the standard module conjecture}\label{subsec langlands classification}

Let $P=MN$ be a parabolic subgroup of $G$ and $\sigma$ a finite length representation of $M$. We let $\Ind_{P}^{G}(\sigma)$ denote the normalized parabolic induction, which is a finite length representation of $G$. If $P'=MN'$ is another parabolic subgroup of $G$ containing $M$, then $\Ind_{P}^{G}(\sigma)$ and $\Ind_{P'}^{G}(\sigma)$ have the same image in the Grothendieck group (\cite[\S 2.9, 2.10]{BZ77}). Thus, we let $\Ind_{M}^{G}(\sigma)$ denote the image of $\Ind_{P}^{G}(\sigma)$ in the Grothendieck group for any choice of $P$. Similarly, for a finite length representation $\pi$ of $G$, we let $\Jac_{P}(\pi)$ denote the normalized Jacquet module of $\pi$, which is a finite length representation of $M$, and let $\Jac_{M}(\pi)$ denote its image in the Grothendieck group.

The Langlands classification for $\Pi(G)$ (\cite[Theorem 3.5]{Kon03}) and for $\Phi(G)$ (\cite{SZ18}) gives canonical bijections
\begin{align*}
    \pi \leftrightarrow (P,\pi_{temp},\nu),\ \ \phi \leftrightarrow (P,\phi_{temp},\nu),
\end{align*}
where
\begin{itemize}
    \item [$\oldbullet$] $\pi \in \Pi(G)$ and $\phi \in \Phi(G)$;
    \item [$\oldbullet$] $P$ is a standard parabolic subgroup of $G$ with standard Levi subgroup $M$;
    \item [$\oldbullet$] $\nu$ is an unramified character of $M$, corresponding to an element of $X_{\ast}(M)\otimes\mathbb{R}$, in the positive Weyl chamber determined by $P$;
    \item [$\oldbullet$] $\pi_{temp}\in \Pi(M)$ and $\phi_{temp}\in \Phi(M)$ are both tempered.
\end{itemize}
Under this bijection, the $L$-parameter $\phi$ is given by
\[ \phi = \left({}^{L}M\hookrightarrow{}^{L}G\right)\circ(\phi_{temp}\otimes\nu),\]
where $\nu$ is regarded as an unramified central cocharacter of $\widehat{M}$, so that $\phi_{temp}\otimes\nu$ is the corresponding unramified twist of $\phi_{temp}$.

On the representation side, $\pi$ is the unique irreducible quotient of $\Ind_{P}^{G}(\pi_{temp}\otimes \nu)$. We call this induced representation the \emph{standard module} of $\pi$ and denote it by $M(\pi)$. If $\pi$ is generic, then $\pi_{temp}$ is also generic.

We shall also need the subrepresentation version of the Langlands classification. Let $(P,\pi_{temp},\nu)$ be a standard triple with Langlands quotient $\pi$. There is an intertwining operator
\[ j(\nu): \Ind_{P}^{G}\pi_{temp}\otimes \nu \longrightarrow \Ind_{\ol{P}}^{G}\pi_{temp}\otimes \nu,\]
where $\ol{P}$ is the opposite parabolic subgroup, such that $\pi \cong \im(j(\nu))$. Thus, we obtain an injection
\begin{align}\label{eq langlands sub}
    \pi \hookrightarrow \Ind_{\ol{P}}^{G}\pi_{temp}\otimes \nu \cong \Ind_{w(\ol{P})}^{G} w(\pi_{temp})\otimes w(\nu),
\end{align}
where $w$ is the Weyl element that sends $\ol{P}$ to a standard parabolic subgroup $w(\ol{P})$. The triple $(w(\ol{P}), w(\pi_{temp}), w(\nu))$ then satisfies all the conditions of a standard triple, except that $w(\nu)$ is negative with respect to $w(\ol{P})$. We let $M_{sub}(\pi)$ denote the right-hand side of \eqref{eq langlands sub}, the standard module of $\pi$ in the subrepresentation setting. By \cite[\S 2.9, 2.10]{BZ77}, the modules $M(\pi)$ and $M_{sub}(\pi)$ have the same image in the Grothendieck group.

Finally, we recall the generalized injectivity conjecture of Casselman and Shahidi, now a theorem.

\begin{theorem}[Generalized injectivity conjecture, {\cite[Conjecture 1.2]{CS98}}; {\cite{Han10, Han20}}, {\cite{Dij22}}, {\cite[Theorem D(b)]{Sol23}}]\label{thm generalized injectivity}
    Let $G$ be quasi-split and let $\tau = \Ind_{P}^{G}\pi_{temp}\otimes \nu$ be a standard module of $G$ with $\pi_{temp}$ generic. Then any generic irreducible subquotient of $\tau$ is a subrepresentation of $\tau$. That is, if $\pi \in \Pi(G)$ is generic and $\pi \leq \tau$, then $\Hom(\pi,\tau)\neq 0$.
\end{theorem}

\begin{remark}\label{rmk injectivity}
    The generalized injectivity conjecture was established in several stages. Casselman and Shahidi (\cite{CS98}) proved it when the inducing datum is supercuspidal. Hanzer (\cite{Han10, Han20}) proved it for the quasi-split classical groups. Dijols (\cite{Dij22}) then proved it whenever $G$ has no simple factor of exceptional type, and for many Bernstein blocks of the remaining groups. Finally, Solleveld (\cite{Sol23}) gave a geometric proof for an arbitrary quasi-split $G$, by translating the question into one about the affine Hecke algebra of a generic Bernstein block.
\end{remark}

\subsection{The Aubert-Zelevinsky involution}\label{subsec AZ}

Let $\mathcal{R}(G)$ be the Grothendieck group of smooth representations of finite length of $G$. If $\pi$ is a smooth representation of finite length of $G$, we let $[\pi]$ denote its image in $\mathcal{R}(G)$.

Aubert (\cite{Aub95}) showed that for any $\pi \in \Pi(G)$, there exists $\varepsilon \in \{\pm 1\}$ and an irreducible representation $\widehat{\pi}\in \Pi(G)$ such that
\begin{align}\label{eq AZ}
    [\widehat{\pi}]:= \varepsilon \sum_{P}(-1)^{\dim(A_P)}\left[\Ind_{P}^{G}\circ \Jac_{P}(\pi)\right].
\end{align}
Here the sum is taken over all standard parabolic subgroups $P$ of $G$ and $A_P$ is the maximal split torus of the center of the Levi subgroup of $P$. Moreover, the map $\pi \mapsto \widehat{\pi}$ is an involution on $\Pi(G)$. We call $\widehat{\pi}$ the \emph{Aubert-Zelevinsky involution} of $\pi$, and we let $\AZ$ denote the linear transformation on the Grothendieck group $\mathcal{R}(G)$, and likewise on $\mathcal{R}(M)$ for any Levi subgroup $M$ of $G$, extending $\pi \mapsto \widehat{\pi}$ on irreducible representations.
The involution preserves supercuspidal support. In particular, $\AZ(\pi)=\pi$ when $\pi$ is supercuspidal.

For a standard Levi subgroup $M$ of $G$ with standard parabolic subgroup $P=MN$, let $w_{\ell}$ and $w_{\ell}^{M}$ be the longest elements of the Weyl group of $G$ and of $M$, respectively, with respect to a fixed maximal split torus $A_0$, and set
\[ w_M := w_{\ell}\, w_{\ell}^{M}. \]
Then $w_M M w_M^{-1}$ is a standard Levi subgroup, and $w_M$ sends the opposite parabolic $\ol{P}$ to the standard parabolic subgroup with Levi $w_M M w_M^{-1}$. We shall need the following compatibility of the Aubert-Zelevinsky involution with Jacquet modules.

\begin{proposition}[{\cite[Th{\'e}or{\`e}me 1.7(2)]{Aub95}}]\label{prop AZ jacquet}
    Let $M$ be a standard Levi subgroup of $G$ with standard parabolic subgroup $P=MN$. Then, for any $\pi \in \Pi(G)$, we have an identity in $\mathcal{R}(M)$
    \[ \AZ\circ \Jac_{M}(\pi) = \Ad(w_M)\inv \circ \Jac_{w_M M w_M^{-1}}\left(\AZ(\pi)\right).\]
\end{proposition}

\subsection{The local Langlands correspondence}\label{subsec LLC}

We work with the following version of the local Langlands correspondence, a system of maps for $G$ and all of its standard Levi subgroups.

\begin{conjecture}[Local Langlands correspondence]\label{conj LLC}
    For each $M \in \mathcal{L}(G)$, there is a canonical map
    \begin{align*}
        \LLC_{M}: \Pi(M) &\longrightarrow \Phi(M),\\
        \pi & \longmapsto \phi_{\pi},
    \end{align*}
    with finite fibers. For $\phi \in \Phi(M)$, we write $\Pi_{\phi}:= \LLC_{M}^{-1}(\phi)$ for the associated $L$-packet.
\end{conjecture}

We write $\LLC_{\mathcal{L}(G)}:= \{\LLC_M\}_{M \in \mathcal{L}(G)}$, and we drop the subscript $M$ when it is clear from the context. Following \cite[\S 3]{KLT}, we record the following two desiderata for $\LLC_{\mathcal{L}(G)}$. They are the only ones used in \S \ref{sec reduction}.

\begin{desideratum}\label{des LC IT}
    Let $\LLC_{\mathcal{L}(G)}$ be a system of maps as in Conjecture \ref{conj LLC}.
    \begin{enumerate}
        \item [\emph{\textbf{(LC)}}] Suppose that $\pi \leftrightarrow (P,\pi_{temp},\nu)$ and $\LLC_{G}(\pi)\leftrightarrow (P',\phi_{temp},\nu')$ under the Langlands classification. Then $P=P'$, $\nu=\nu'$ and $\phi_{temp}= \LLC_{M}(\pi_{temp})$.
        \item [\emph{\textbf{(IT)}}] Suppose $\pi \leq \Ind_{P}^{G}\sigma$, where $\sigma$ is a (unitary) tempered representation of $M$, the Levi subgroup of $P$. Then
        \[ \phi_{\pi}= \left({}^{L}M \hookrightarrow {}^{L}G\right)\circ \phi_{\sigma}.\]
    \end{enumerate}
\end{desideratum}

Finally, we recall the infinitesimal parameter attached to an $L$-parameter. For $\phi \in \Phi(G)$, define $\lambda_{\phi}: W_F \to {}^{L}G$ by
\begin{align}\label{eq inf param}
    \lambda_{\phi}(w):= \phi\left(w, \begin{pmatrix} |w|^{1/2} & \\ & |w|^{-1/2}\end{pmatrix}\right),
\end{align}
which we regard as an $L$-parameter of $G$ trivial on $\SL_2(\BC)$.
Fixing an infinitesimal parameter $\lambda$ of $G$, the set
\[ \Phi_{\lambda}(G):= \{ \phi \in \Phi(G)\ | \ \lambda_{\phi}= \lambda\}\]
is finite and is in bijection with the set of orbits of the Vogan variety $V_{\lambda}$ associated to $\lambda$; see \cite[\S 4]{CFMMX22}.  Note that $\phi$ is trivial on $\SL_2(\BC)$ if and only if $\phi = \lambda_{\phi}$, which corresponds to the zero orbit in $V_{\lambda}$.
The variety $V_{\lambda}$ has a unique open orbit, and we say that the corresponding $L$-parameter is \emph{open}. 
Openness can also be read off from the adjoint $L$-function. 

\begin{proposition}[{\cite[Proposition 4.4]{CDFZ24}}]\label{prop open adjoint L}
    Let $\phi$ be an $L$-parameter of $G$. The adjoint $L$-function $L(s,\phi,\Ad)$ is regular at $s=1$ if and only if $\phi$ is open.
\end{proposition}

\section{Reduction to the supercuspidal case}\label{sec reduction}

In this section we prove Theorem \ref{thm reduction}.

\begin{proof}[Proof of Theorem \ref{thm reduction}]
    Let $\pi$ be a generic representation of $G$. Choose a unitary supercuspidal representation $\pi_{sc}$ of a standard Levi subgroup $M$ and an unramified character $\nu$ of $M$ such that
    \[ \pi \leq \Ind_{P}^{G} \pi_{sc}\otimes \nu.\]
    Let $P$ br the standard parabolic subgroup $P$ for $M$. We may assume that $\nu$ is non-negative with respect to $P$ by applying a Weyl action if necessary. 
    Let $P'=M'N'$ be the larger standard parabolic subgroup for which $\nu$ extends to an unramified character $\nu'$ of $M'$ that is positive with respect to $P'$. Then
    \[ \Ind_{P}^{G} \pi_{sc}\otimes \nu = \Ind_{P'}^{G}\left( \Ind_{P\cap M'}^{M'} \pi_{sc}\right)\otimes \nu'.\]
    Let $\pi_{temp}$ be the unique $\mfr{w}_{M'}$-generic summand of the unitary tempered induction $\Ind_{P\cap M'}^{M'}\pi_{sc}$, which exists by Proposition \ref{prop rodier}. Then we have
    \[ \pi \leq \Ind_{P'}^{G} \pi_{temp}\otimes \nu' = M(\pi').\]
    By Theorem \ref{thm generalized injectivity}, we have
    \[ \Hom\left(\pi, M(\pi')\right)\neq 0.\]
    Composing with the injection induced from $\pi_{temp}\hookrightarrow \Ind_{P\cap M'}^{M'}\pi_{sc}$, we obtain
    \[ \Hom\left(\pi, \Ind_{P}^{G} \pi_{sc}\otimes \nu\right)\neq 0.\]
    Now Frobenius reciprocity gives
    \[ \Hom\left( \Jac_{P}(\pi), \pi_{sc}\otimes \nu\right)\neq 0.\]
    This implies that $\AZ(\pi_{sc}\otimes \nu)= \pi_{sc}\otimes \nu$ is an irreducible subquotient of
    \[ \AZ\circ \Jac_{M}(\pi) = \Ad(w_M)\inv \circ \Jac_{\widetilde{M}}\left( \AZ(\pi)\right),\]
    where the equality is Proposition \ref{prop AZ jacquet} and $\widetilde{M}= w_M M w_M^{-1}$ is a standard Levi subgroup, with standard parabolic subgroup $\widetilde{P}$.

    As $\pi_{sc}\otimes \nu$ is supercuspidal, it appears as a quotient of the finite length representation $\Ad(w_M)\inv\circ \Jac_{\widetilde{P}}(\widehat{\pi})$, of which it is a subquotient (\cite[Proposition 3.30]{BZ76}, \cite[Lemma VI.3.6]{Ren26}). Thus,
    \begin{align*}
        0 \neq \Hom\left( \Ad(w_M)\inv \circ \Jac_{\widetilde{P}}(\widehat{\pi}),\ \pi_{sc}\otimes\nu\right)
        \cong \Hom\left( \Jac_{\widetilde{P}}(\widehat{\pi}),\ \Ad(w_M)\cdot (\pi_{sc}\otimes \nu)\right).
    \end{align*}
    By Frobenius reciprocity, we get
    \[ 0 \neq \Hom\left( \widehat{\pi},\ \Ind_{\widetilde{P}}^{G} \pi_{sc}^{w_M}\otimes w_M(\nu)\right),\]
    where $\pi_{sc}^{w_M}\otimes w_M(\nu)= \Ad(w_M)(\pi_{sc}\otimes \nu)$ is a supercuspidal representation of $\widetilde{M}$.
    Since $w_{\ell}^{M}$ fixes the split component $A_M$ pointwise, we have $w_M(\nu)= w_{\ell}\, w_{\ell}^{M}(\nu)= w_{\ell}(\nu)$. As $w_{\ell}$ is the longest Weyl element and $\nu$ is non-negative with respect to $P$, the character $w_M(\nu)$ is non-positive with respect to $\widetilde{P}$. Let $P''=M''N''$ be the standard parabolic subgroup containing $\widetilde{P}$ for which $w_M(\nu)$ is negative with respect to $P''$, and choose a tempered representation
    \[ \pi_{temp}'' \leq \Ind_{\widetilde{P}\cap M''}^{M''} \pi_{sc}^{w_M}\]
    such that
    \[ 0 \neq \Hom\left( \widehat{\pi},\ \Ind_{P''}^{G} \pi_{temp}''\otimes w_M(\nu)\right).\]
    The right hand side is the subrepresentation standard module $M_{sub}(\widehat{\pi})$. By the argument around \eqref{eq langlands sub}, the Langlands classification assigns to $\widehat{\pi}$ the standard triple
    \[ \left( w_{M''}(\ol{P''}),\ w_{M''}(\pi_{temp}''),\ w_{M''}w_M(\nu)\right). \]
    Thus \textbf{(LC)} gives
    \begin{align}\label{eq: reduction 2}
         \LLC(\widehat{\pi}) = \left( {}^{L}\bigl(w_{M''}(M'')\bigr) \hookrightarrow {}^{L}G\right)\circ \LLC\!\left(w_{M''}(\pi_{temp}'')\right)\otimes w_{M''}w_M(\nu).
    \end{align}

    Now recall that $\pi_{temp}'' \leq \Ind_{\widetilde{P}\cap M''}^{M''}\pi_{sc}^{w_M}$. Conjugating by $w_{M''}$, we obtain
    \[ w_{M''}(\pi_{temp}'') \leq \Ind_{w_{M''}(\widetilde{P}\cap M'')}^{w_{M''}(M'')} \sigma, \qquad \sigma := w_{M''}\!\left(\pi_{sc}^{w_M}\right), \]
    so that $w_{M''}(\pi_{temp}'')$ is a subquotient of a parabolic induction of the supercuspidal representation $\sigma$ of $w_{M''}(\widetilde{M})$. By \textbf{(IT)},
    \begin{align}\label{eq: reduction 1}
        \LLC\!\left(w_{M''}(\pi_{temp}'')\right) = \left( {}^{L}\bigl(w_{M''}(\widetilde{M})\bigr) \hookrightarrow {}^{L}\bigl(w_{M''}(M'')\bigr)\right)\circ \LLC(\sigma).
    \end{align}

    Finally, $\sigma$ is a Weyl conjugate of the generic supercuspidal representation $\pi_{sc}$, so it is generic with respect to some Whittaker datum. By the hypothesis of the theorem, Conjecture \ref{conj main} holds for $\sigma$. We conclude that 
    \[ \LLC(\widehat{\sigma})= \LLC(\sigma)\]
     is trivial on $\SL_2(\BC)$. Then it follows from \eqref{eq: reduction 1} and \eqref{eq: reduction 2} that $\LLC\!\left(w_{M''}(\pi_{temp}'')\right)$ and $\LLC(\widehat{\pi})$ are both  trivial on $\SL_2(\BC)$ as well. This proves Conjecture \ref{conj main} for $\pi$ and completes the proof of the theorem.
\end{proof}

\begin{remark}\label{rmk weyl}
    Suppose in addition that $\LLC$ is compatible with the Weyl group action. Then \eqref{eq: reduction 1} and \eqref{eq: reduction 2} in the proof above show that $\LLC(\widehat{\pi})$ is conjugate to
    \[ \left( {}^{L}M'' \hookrightarrow {}^{L}G\right)\circ \left(\LLC(\pi_{temp}'')\otimes w_M(\nu)\right).\]
\end{remark}

\section{The supercuspidal case}\label{sec supercuspidal}

In this section we give two sets of assumptions on $\LLC$, each of which implies Conjecture \ref{conj main} for generic supercuspidal representations of the Levi subgroups of $G$. Together with Theorem \ref{thm reduction}, this proves Theorem \ref{thm intro main}.

\subsection{Compatibility with Kaletha's local Langlands correspondence}\label{subsec kaletha}

In \cite{KLT}, we prove the following theorem.

\begin{theorem}[\cite{KLT}]\label{thm KLT}
    If $p \gg 0$, then every generic supercuspidal representation is non-singular in the sense of \cite{Kal21}.
\end{theorem}

In \cite{Kal21}, Kaletha defined an $\LLC$:
\begin{align*}
    \LLC^{\Ka}: \left\{\begin{matrix} \text{Non-singular supercuspidal}\\ \text{representations of } G\end{matrix}\right\}
    \longrightarrow
    \left\{\begin{matrix} \text{discrete } L\text{-parameters of}\\ G \text{ trivial on } \SL_2(\BC)\end{matrix}\right\}.
\end{align*}
 Thus, if we require that $\LLC(\pi)= \LLC^{\Ka}(\pi)$ whenever $\pi$ is non-singular, then Theorem \ref{thm KLT} shows that Conjecture \ref{conj main} holds for all generic supercuspidal representations of the Levi subgroups of $G$ when $p \gg 0$.

\subsection{Conjectures of Gross-Prasad-Rallis and Haines}\label{subsec GPR Haines}

 In this subsection we give a second argument for the supercuspidal case without assumptions on $p$. Instead, we assume three more properties of $\LLC$ that are expected to hold in general. The first two describe the generic members of a single $L$-packet, and the third describes how infinitesimal parameters behave under parabolic induction.

The first is a refinement of the Shahidi conjecture, due to Gross-Prasad and Rallis. It says that the existence of a generic member in an $L$-packet is detected by the adjoint $L$-function of its parameter.

\begin{conjecture}[Gross-Prasad-Rallis, {\cite[Conjecture 2.6]{GP92}}]\label{conj GPR}
    Let $G$ be quasi-split and fix a Whittaker datum $\mfr{w}$. Then $L(s,\phi,\Ad)$ is regular at $s=1$ if and only if $\Pi_{\phi}$ contains a $\mfr{w}$-generic member.
\end{conjecture}

The second says that such a member, when it exists, is unique. 
\begin{conjecture}\label{conj uniqueness}
    Let $G$ be quasi-split and fix a Whittaker datum $\mfr{w}$. Then for any tempered $\phi \in \Phi(G)$, the $L$-packet $\Pi_{\phi}$ contains at most one $\mfr{w}$-generic member.
\end{conjecture}

\begin{remark}\label{rmk GPR uniqueness}
    The existence assertion in Conjecture \ref{conj GPR} goes back to Shahidi (\cite[Conjecture 9.4]{Sha90}). Conjecture \ref{conj uniqueness} is the stronger form proved by Shelstad over the archimedean fields (\cite{She08}). The non-archimedean analogue is expected in general. It is proved for Arthur's local Langlands correspondence for symplectic and split odd special orthogonal groups by Atobe (\cite{Ato17}).

    By Proposition \ref{prop open adjoint L}, the condition that $L(s,\phi,\Ad)$ be regular at $s=1$ may equivalently be replaced by the condition that $\phi$ be open, and we shall use the two formulations interchangeably.
\end{remark}

We extend the uniqueness of generic member to $L$-packet with open parameter that is not necessarily tempered by \textbf{(LC)}.

\begin{lemma}\label{lem uniqueness open}
    Assume that $\LLC_{\mathcal{L}(G)}$ satisfies \textbf{(LC)} and that Conjecture \ref{conj uniqueness} holds for every standard Levi subgroup of $G$. Then for any open $\phi \in \Phi(G)$, the $L$-packet $\Pi_{\phi}$ contains at most one $\mfr{w}$-generic member.
\end{lemma}
\begin{proof}
    Write $\phi \leftrightarrow (P=MN, \phi_{temp},\nu)$ under the Langlands classification, so that $\phi_{temp}\in \Phi(M)$ is tempered. By \textbf{(LC)}, the $L$-packet $\Pi_{\phi}$ is the collection of Langlands quotients $\pi$ of $\Ind_{P}^{G}\pi_{temp}\otimes \nu$, where $\pi_{temp}$ ranges over $\Pi_{\phi_{temp}}$. Moreover, $\pi \mapsto \pi_{temp}$ is a bijection between $\Pi_{\phi}$ and $\Pi_{\phi_{temp}}$.

    If $\pi$ is $\mfr{w}$-generic, then $\Ind_{P}^{G}\pi_{temp}\otimes \nu$ has a $\mfr{w}$-generic subquotient, so by the defining property of $\mfr{w}_M$ in \S \ref{subsec whittaker} the representation $\pi_{temp}\otimes \nu$, and hence $\pi_{temp}$, is $\mfr{w}_M$-generic. As $\phi_{temp}$ is tempered, Conjecture \ref{conj uniqueness} for $M$ says $\Pi_{\phi_{temp}}$ has at most one $\mfr{w}_M$-generic member. Through the bijection $\pi \mapsto \pi_{temp}$, it follows that $\Pi_{\phi}$ has at most one $\mfr{w}$-generic member.
\end{proof}

Finally, we recall the compatibility of the local Langlands correspondence with supercuspidal support proposed by Haines. Recall the infinitesimal parameter $\lambda_{\phi}$ from \eqref{eq inf param}.

\begin{conjecture}[Haines, {\cite[Conjecture 5.2.2]{Hai14}}]\label{conj Haines}
    Let $Q$ be a parabolic subgroup of $G$ with Levi subgroup $M$ and let $\tau \in \Pi(M)$. If $\pi \leq \Ind_{Q}^{G}\tau$, then
    \[ \lambda_{\LLC(\pi)} = \left({}^{L}M \hookrightarrow {}^{L}G\right)\circ \lambda_{\LLC(\tau)}.\]
\end{conjecture}

We now prove Conjecture \ref{conj main} assuming the above three expected properties of $\LLC$.

\begin{theorem}\label{thm supercuspidal}
    Assume that $\LLC_{\mathcal{L}(G)}$ satisfies \textbf{(LC)} and Conjectures \ref{conj GPR}, \ref{conj uniqueness} and \ref{conj Haines}. Then Conjecture \ref{conj main} holds for supercuspidal representations.
\end{theorem}

\begin{proof}
    We argue by contradiction. Suppose that $\pi$ is generic supercuspidal but $\phi := \LLC(\pi)$ is not trivial on $\SL_2(\BC)$. Write
    \[ \lambda := \lambda_{\phi} \longleftrightarrow (P= MN, \lambda_{temp},\nu)\]
    under the Langlands classification. Since $\lambda_{\phi}\neq \phi$ by assumption, $\lambda$ is not open, and so $M$ must be proper.

    By the Gross-Prasad-Rallis conjecture, there is a $\mfr{w}_M$-generic $\pi_{sc}\in \Pi_{\lambda_{temp}}$. Hence the induced representation $\Ind_{P}^{G}\pi_{sc}\otimes \nu$ contains a $\mfr{w}$-generic subquotient $\pi_{gen}$, which is not supercuspidal. The Gross-Prasad-Rallis conjecture also says that $L(s,\LLC(\pi_{gen}),\Ad)$ is regular at $s=1$, so $\LLC(\pi_{gen})$ is open by Proposition \ref{prop open adjoint L}. On the other hand, Haines' conjecture gives
    \[ \lambda_{\LLC(\pi_{gen})} = \left({}^{L}M \hookrightarrow {}^{L}G\right)\circ (\lambda_{temp}\otimes \nu) = \lambda.\]
    These two facts together imply that $\LLC(\pi_{gen})= \phi$, the unique open $L$-parameter in $\Phi_{\lambda}(G)$.
 However, then $\pi$ and $\pi_{gen}$ are two $\mfr{w}$-generic members of the $L$-packet $\Pi_{\phi}$, and they are distinct since $\pi$ is supercuspidal and $\pi_{gen}$ is not. As $\phi$ is open, this contradicts Lemma \ref{lem uniqueness open}. We conclude that $\phi$ must be trivial on $\SL_2(\BC)$.
\end{proof}

This completes the proof of Theorem \ref{thm intro main}. Indeed, \textbf{(LC)} and \textbf{(IT)} let Theorem \ref{thm reduction} reduce Conjecture \ref{conj main} to the generic supercuspidal representations of the Levi subgroups of $G$, and this case is handled by \S \ref{subsec kaletha} under hypothesis (a) and by Theorem \ref{thm supercuspidal} under hypothesis (b).


\begin{thebibliography}{CFMMX22}



\bibitem[Ato17]{Ato17} H. Atobe, \emph{On the uniqueness of generic representations in an $L$-packet}, Int. Math. Res. Not. IMRN \textbf{2017}, no. 23, 7051--7068.

\bibitem[AC26]{AC26} H. Atobe and D. Ciubotaru, \emph{Endoscopic transfer and the wavefront upper bound conjecture}, (2026), arXiv:2602.22504.

\bibitem[Aub95]{Aub95} A.-M. Aubert, \emph{Dualit{\'e} dans le groupe de Grothendieck de la cat{\'e}gorie des repr{\'e}sentations lisses de longueur finie d'un groupe r{\'e}ductif $p$-adique}, Trans. Amer. Math. Soc. \textbf{347} (1995), 2179--2189.

\bibitem[BV85]{BV85} D. Barbasch and D. Vogan, \emph{Unipotent representations of complex semisimple groups}, Ann. of Math. (2) \textbf{121} (1985), no. 1, 41--110.

\bibitem[BZ76]{BZ76} I. N. Bernstein and A. V. Zelevinsky, \emph{Representations of the group $GL(n,F)$ where $F$ is a non-Archimedean local field}, Russian Math. Surveys \textbf{31} (1976), no. 3, 1--68.

\bibitem[BZ77]{BZ77} I. N. Bernstein and A. V. Zelevinsky, \emph{Induced representations of reductive $p$-adic groups I}, Ann. Sci. {\'E}cole Norm. Sup. (4) \textbf{10} (1977), no. 4, 441--472.

\bibitem[CS80]{CS80} W. Casselman and J. Shalika, \emph{The unramified principal series of $p$-adic groups. II. The Whittaker function}, Compositio Math. \textbf{41} (1980), no. 2, 207--231.

\bibitem[CS98]{CS98} W. Casselman and F. Shahidi, \emph{On irreducibility of standard modules for generic representations}, Ann. Sci. {\'E}cole Norm. Sup. (4) \textbf{31} (1998), no. 4, 561--589.


\bibitem[CK25]{CK25} D. Ciubotaru and J.-L. Kim, \emph{The wavefront set: bounds for the Langlands parameter}, Math. Ann. \textbf{393} (2025), no. 2, 1827--1861.

\bibitem[CDFZ24]{CDFZ24} C. Cunningham, S. Dijols, A. Fiori, and Q. Zhang, \emph{Generic representations, open parameters and ABV-packets for $p$-adic groups}, (2024), arXiv:2404.07463.


\bibitem[CFMMX22]{CFMMX22} C. Cunningham, A. Fiori, A. Moussaoui, J. Mracek, and B. Xu, \emph{Arthur packets for $p$-adic groups by way of microlocal vanishing cycles of perverse sheaves, with examples}, Mem. Amer. Math. Soc. \textbf{276} (2022), no. 1353.




\bibitem[Dij22]{Dij22} S. Dijols, \emph{The generalized injectivity conjecture}, Bull. Soc. Math. France \textbf{150} (2022), no. 2, 251--345.

\bibitem[GP92]{GP92} B. H. Gross and D. Prasad, \emph{On the decomposition of a representation of $\SO_n$ when restricted to $\SO_{n-1}$}, Canad. J. Math. \textbf{44} (1992), no. 5, 974--1002.

\bibitem[Hai14]{Hai14} T. Haines, \emph{The stable Bernstein center and test functions for Shimura varieties}, in: Automorphic forms and Galois representations, Vol. 2, London Math. Soc. Lecture Note Ser. \textbf{415}, Cambridge Univ. Press, 2014, 118--186.

\bibitem[Han10]{Han10} M. Hanzer, \emph{The generalized injectivity conjecture for the classical $p$-adic groups}, Int. Math. Res. Not. IMRN \textbf{2010}, no. 2, 195--237.

\bibitem[Han20]{Han20} M. Hanzer, \emph{Generalized injectivity conjecture for classical $p$-adic groups II}, J. Pure Appl. Algebra \textbf{224} (2020), no. 1, 149--168.

\bibitem[HC78]{HC78} Harish-Chandra, \emph{Admissible invariant distributions on reductive $p$-adic groups}, in: Lie theories and their applications, Queen's Papers in Pure Appl. Math. \textbf{48}, Queen's Univ., Kingston, Ont., 1978, 281--347.

\bibitem[HLLS24]{HLLS24} A. Hazeltine, B. Liu, C.-H. Lo, and F. Shahidi, \emph{On the upper bound of wavefront sets of representations of $p$-adic groups}, (2024), arXiv:2403.11976.

\bibitem[Jia14]{Jia14} D. Jiang, \emph{Automorphic integral transforms for classical groups I: endoscopy correspondences}, in: Automorphic forms: $L$-functions and related geometry, Contemp. Math. \textbf{614}, Amer. Math. Soc., 2014, 179--242.

\bibitem[Kal21]{Kal21} T. Kaletha, \emph{Supercuspidal $L$-packets}, (2021), arXiv:1912.03274.

\bibitem[KLT]{KLT} T. Kaletha, C.-H. Lo and C.-C. Tsai, \emph{On desiderata and uniqueness of the local Langlands correspondence}, in preparation.

\bibitem[Kon03]{Kon03} T. Konno, \emph{A note on the Langlands classification and irreducibility of induced representations of $p$-adic groups}, Kyushu J. Math. \textbf{57} (2003), no. 2, 383--409.

\bibitem[LLS24]{LLS24} B. Liu, C.-H. Lo, and F. Shahidi, \emph{On anti-tempered local Arthur packets and a lemma of Arthur}, (2024), arXiv:2405.17407.

\bibitem[LS23]{LS23} B. Liu and F. Shahidi, \emph{Jiang's conjecture on the wavefront sets of local Arthur packets}, submitted, 2023.

\bibitem[MW87]{MW87} C. M{\oe}glin and J.-L. Waldspurger, \emph{Mod{\`e}les de Whittaker d{\'e}g{\'e}n{\'e}r{\'e}s pour des groupes $p$-adiques}, Math. Z. \textbf{196} (1987), 427--452.

\bibitem[Ren26]{Ren26} D. Renard, \emph{Representations of $p$-adic reductive groups}, (2026), arXiv:2607.08308.

\bibitem[Rod73]{Rod73} F. Rodier, \emph{Whittaker models for admissible representations of reductive $p$-adic split groups}, in: Harmonic analysis on homogeneous spaces, Proc. Sympos. Pure Math. \textbf{26}, Amer. Math. Soc., 1973, 425--430.

\bibitem[Sha90]{Sha90} F. Shahidi, \emph{A proof of Langlands' conjecture on Plancherel measures; complementary series for $p$-adic groups}, Ann. of Math. (2) \textbf{132} (1990), no. 2, 273--330.

\bibitem[She08]{She08} D. Shelstad, \emph{Tempered endoscopy for real groups III: inversion of transfer and $L$-packet structure}, Represent. Theory \textbf{12} (2008), 369--402.


\bibitem[SZ18]{SZ18} A. J. Silberger and E.-W. Zink, \emph{Langlands classification for $L$-parameters}, J. Algebra \textbf{511} (2018), 299--357.

\bibitem[Sol23]{Sol23} M. Solleveld, \emph{On submodules of standard modules}, (2023), arXiv:2309.10401.



\bibitem[Var14]{Var14} S. Varma, \emph{On a result of M{\oe}glin and Waldspurger in residual characteristic $2$}, Math. Z. \textbf{277} (2014), no. 3--4, 1027--1048.

\end{thebibliography}
\end{document}